\documentclass[a4paper, 11pt, reqno, onesided]{amsart}

\usepackage{amssymb,amsmath}
\usepackage[dvips]{graphics}
\usepackage[dvips]{color}
\usepackage{graphicx}
\usepackage[english]{babel}
\usepackage[latin1]{inputenc}
\usepackage{comment}
\usepackage{mathrsfs}

\newtheorem{theorem}{Theorem}[section]
\newtheorem{lemma}[theorem]{Lemma}

\newtheorem{corollary}[theorem]{Corollary}

\theoremstyle{definition}
\newtheorem{definition}[theorem]{Definition}

\newtheorem{remark}[theorem]{Remark}

\numberwithin{equation}{section}

\begin{document}

\title{The equivalence between pointwise Hardy inequalities and uniform fatness}

\author{Riikka Korte \and Juha Lehrb\"ack \and Heli Tuominen} 
\renewcommand{\datename}{\textit{}}
\thanks{The second and the third author were supported partially by
the Academy of Finland, grant no.\ 120972}

\newcommand\rn{\mathbb R^n}
\newcommand\re{\mathbb R}
\newcommand\R{\mathbb R}
\newcommand\N{\mathbb N}
\newcommand\lp{L^p(X)}
\newcommand\grad{\nabla}
\newcommand\bd{\partial}
\newcommand\eps{\varepsilon}
\newcommand\Hc{\mathcal H}
\newcommand\Hh{\mathscr H}
\newcommand\Fc{\mathcal F}
\newcommand\M{\operatorname{\mathcal M}}
\newcommand\diam{\operatorname{diam}}
\newcommand\rad{\operatorname{rad}}
\newcommand\dist{\operatorname{d}}
\newcommand\len{\operatorname{\ell}}
\newcommand\ol{\overline}
\newcommand\sub{\subset}
\newcommand\vphi{\varphi}
\newcommand{\dom}{\delta_{\Omega}}

\providecommand{\pcap}[1]{\operatorname{cap}_p({#1})}
\providecommand{\hcont}[1]{\operatorname{\Hh}^h_{\infty}({#1})}
\providecommand{\sob}[1]{W^{1,p}({#1})}
\providecommand{\sobn}[1]{N^{1,p}({#1})}
\providecommand{\sobno}[1]{N_0^{1,p}({#1})}
\providecommand{\smooth}[1]{C_0^{\infty}({#1})}
\providecommand{\dic}[1]{\operatorname{d}_{\Omega}({#1})}
\providecommand{\ch}[1]{\text{\raise 2pt \hbox{$\chi$}\kern-0.2pt}_{#1}}
\providecommand{\vint}[1]{\mathchoice
          {\mathop{\vrule width 5pt height 3 pt depth -2.5pt
                  \kern -8pt \intop}\nolimits_{\kern -5pt{#1}}}%
          {\mathop{\vrule width 5pt height 3 pt depth -2.6pt
                  \kern -6pt \intop}\nolimits_{\kern -3pt{#1}}}%
          {\mathop{\vrule width 5pt height 3 pt depth -2.6pt
                  \kern -6pt \intop}\nolimits_{\kern -3pt{#1}}}%
          {\mathop{\vrule width 5pt height 3 pt depth -2.6pt
                  \kern -6pt \intop}\nolimits_{\kern -3pt{#1}}}}

\begin{abstract}
We prove an equivalence result between 
the validity of a pointwise Hardy inequality in a domain
and uniform capacity density of the complement. 
This result is new even in Euclidean spaces, but 
our methods apply in general metric spaces as well.
We also present a new transparent proof for the fact that uniform
capacity density implies the classical integral version of the Hardy inequality
in the setting of metric spaces. 
In addition, we consider the relations between the above concepts and 
certain Hausdorff content conditions. 
\end{abstract}

\maketitle

\section{Introduction}\label{intro}

Let $\Omega\subsetneq\R^n$ be a domain and let $u\in C_0^\infty(\Omega)$.
The inequality
\begin{equation}\label{pointwise hardy}
|u(x)|\le C\dist(x,\bd\Omega)
\bigl(\,\M_{2\dist(x,\bd\Omega)}|\grad u|^p(x)\bigr)^{1/p},\quad x\in\Omega,
\end{equation}
where $\M_R$ is the restricted Hardy--Littlewood maximal operator
and $1\leq p<\infty$,
can be viewed as a pointwise variant of
the classical $p$-Hardy inequality 
\begin{equation}\label{integral hardy}
\int_\Omega\frac{|u(x)|^p}{\dist(x,\bd\Omega)^p}\,dx
\le C\int_\Omega|\grad u(x)|^p\,dx.
\end{equation} 
We say that the domain $\Omega$ admits the pointwise $p$-Hardy inequality,
if there exists a constant $C>0$ such that inequality \eqref{pointwise hardy} holds
for all $u\in C_0^\infty(\Omega)$ at every $x\in\Omega$.
As our main result, we prove the following characterization for such domains.
Recall that uniform $p$-fatness is a capacity density condition; 
the exact definition is given in Section \ref{sect: preliminaries}.

\begin{theorem}\label{thm: main} Let $1\leq p<\infty$.
A domain $\Omega\sub\rn$ admits the pointwise $p$-Hardy inequality 
if and only if the complement of $\Omega$ is
uniformly $p$-fat.
\end{theorem}

The origins of Hardy inequalities lie in the one-dimensional considerations by Hardy,
see \cite[\textsection 330]{HLP} and the references therein.
In $\R^n$, for $n\geq 2$, Hardy-type inequalities first appeared in the paper  
of Ne\v{c}as \cite{Ne} in the context of Lipschitz domains.
However, 
it has been well-known since the works of Ancona \cite{An} ($p=2$), Lewis \cite{L}, and
Wannebo \cite{W}, that the regularity of the boundary is not essential for Hardy inequalities. Indeed, uniform $p$-fatness of the complement suffices for a domain 
 to admit the integral $p$-Hardy inequality \eqref{integral hardy}.
Uniform $n$-fatness of the complement is also necessary for the
$n$-Hardy inequality, see \cite{An} and \cite{L},
but this is not true for $p<n$.

Pointwise Hardy inequalities were introduced by Haj\l asz \cite{Hj3} and
Kinnunen and Martio \cite{KiMa2}. 
In these works it was shown that uniform $p$-fatness of the complement
guarantees that 
the domain admits even the pointwise $p$-Hardy inequality;
this is the sufficiency part of Theorem \ref{thm: main}.
 
Using the boundedness of the Hardy--Littlewood
maximal operator it is easy to see that 
a pointwise $q$-Hardy inequality for some $q<p$
implies the $p$-Hardy inequality \eqref{integral hardy}. 
This method does not work if we start with a pointwise $p$-Hardy inequality, 
as only weak type estimates are available when the exponent is not allowed to increase.
Indeed, it has been an open question since the first appearance of pointwise Hardy inequalities
whether the pointwise $p$-Hardy inequality implies the integral $p$-Hardy inequality with
the same exponent.

Now, by a remarkable result of Lewis \cite{L}, uniform $p$-fatness has
the following self-improvement property: If $1<p<\infty$ and a set
$E\sub\R^n$ is uniformly $p$-fat, then $E$ is also uniformly $q$-fat
for some $1<q<p$. Thus Theorem \ref{thm: main} has the striking consequence that
pointwise $p$-Hardy inequalities, for $1<p<\infty$,  
enjoy this same
property. In particular, we obtain a positive answer to the above
question: 

\begin{corollary}\label{coro: pw implies int euk}
Let $1<p<\infty$. If a domain $\Omega\sub\R^n$ admits the pointwise $p$-Hardy inequality,
then $\Omega$ admits the integral $p$-Hardy inequality.
\end{corollary}

In fact, by using the approach of Wannebo, we obtain for Corollary 
\ref{coro: pw implies int euk}
another proof in which we avoid the use of the rather deep
self-improvement of uniform fatness; see Section \ref{sect: from fat to int}. 
In addition, we establish a further equivalence between the conditions of
Theorem \ref{thm: main} and certain Poincar\'e type boundary conditions,
see Theorem \ref{thm: eqv}. Notice also the inclusion of the case $p=1$ in
Theorem \ref{thm: main}. On the contrary, the usual $1$-Hardy inequality
does not hold even in smooth domains.

We remark that it was recently shown in \cite{Le1} that a domain $\Omega\sub\rn$
admits a pointwise $q$-Hardy inequality for some $1<q<p$ if and only if
the complement of $\Omega$ is uniformly $p$-fat
(note here the difference between our terminology and that of \cite{Le1}). 
This result is nevertheless significantly weaker than Theorem \ref{thm: main}, 
as the crucial end-point $q=p$ is not reached.

The second purpose of this paper is to generalize parts of the
existing theory of Euclidean Hardy inequalities to the setting of
metric measure spaces. As a part of this scheme we also state and prove 
Theorem \ref{thm: main} in this more general setting. 
The relevant parts of the analysis in metric spaces, as well as the
exact formulations of our main results, can be found in
Section \ref{sect: preliminaries}. In Section \ref{sect: from fat to pw}, we
prove that uniform $p$-fatness of the complement implies the pointwise
$p$-Hardy inequality also in metric spaces. The necessity part of Theorem 
\ref{thm: main} is then obtained in Section \ref{sect: from pw to fat}. 
Section \ref{sect: from fat to int} contains a transparent proof
for the fact that uniform $p$-fatness of the complement (and thus also
the pointwise $p$-Hardy inequality) is sufficient for $\Omega$ to
admit the usual integral version of the $p$-Hardy inequality. 
Finally, in Section \ref{sect: contents}, we give further
generalizations of the results from \cite{Le1} to metric spaces by linking
pointwise Hardy inequalities and uniform fatness to certain Hausdorff
content density conditions. In the special case of
Carnot--Carath\'eodory spaces similar generalizations were recently obtained
in \cite{DGP}. Different aspects of
Hardy inequalities in the metric setting have also been studied  
in \cite{BMS}, \cite{KiKiMa}, \cite{KoSh}, and \cite{KZ}.

\section{Preliminaries and the main results}\label{sect: preliminaries}

\subsection{Metric spaces}
We assume that $X=(X,\dist,\mu)$ is a complete metric measure space
equipped with a metric $\dist$ and a Borel regular outer measure $\mu$
such that $0<\mu(B)<\infty$ for all balls $B=B(x,r)=\{y\in X:\dist(y,x)<r\}$.
For $0<t<\infty$, we write $tB=B(x,tr)$, and $\ol B$ is the corresponding
closed ball.
We assume that $\mu$ is {\it doubling}, which means that there is
a constant $c_D\ge1$, called {\it the doubling constant of $\mu$}, such that 
\[
\mu(2B)\le c_D\,\mu(B)
\]
for all balls $B$ of $X$. Note that the doubling condition together
with completeness implies that the space is proper, that is, closed
balls of $X$ are compact. 

The doubling condition gives an upper bound for the dimension of $X$.
By this we mean that there is a constant $C=C(c_D)>0$ such that, for $s=\log_2c_D$,
\begin{equation}\label{doubling dimension}
\frac{\mu(B(y,r))}{\mu(B(x,R))}\ge C\Bigl(\frac rR\Bigr)^s
\end{equation}
whenever $0<r\le R<\diam X$ and $y\in B(x,R)$. Inequality
\eqref{doubling dimension} may hold for some smaller exponents
than $\log_2c_D$, too. In such cases we let $s$ denote the
infimum of the exponents for which \eqref{doubling dimension} 
holds and say that $s$ is {\it the doubling dimension} of $X$.

When $\Omega\sub\rn$, we obtain, by the density of smooth functions in
the Sobolev space $W_0^{1,p}(\Omega)$, that the Hardy inequality
\eqref{integral hardy} holds for all $u\in W^{1,p}_0(\Omega)$
if it holds for all smooth functions $\vphi\in C_0^\infty(\Omega)$.
General metric spaces lack the notion of smooth functions, 
but there exists a natural counterpart of Sobolev spaces, defined by 
Shanmugalingam \cite{Sh2} and based on the use of {\it upper gradients}.
We say that a Borel function $g\ge0$ is an upper gradient of a function $u$
on an open set $\Omega\subset X$, if for all curves
$\gamma$ joining points $x$ and $y$ in $\Omega$ we have
\begin{equation}\label{ug}
|u(x)-u(y)|\le\int_{\gamma}g\,ds,
\end{equation}
whenever both $u(x)$ and $u(y)$ are finite, and
$\int_{\gamma}g\,ds=\infty$ otherwise. 
By a {\it curve} we mean a nonconstant, rectifiable, continuous
mapping from a compact interval to $X$. 

If $g\ge0$ is a measurable function and \eqref{ug} only fails for a
curve family with zero $p$-modulus, then $g$ is {\it a $p$-weak upper gradient} 
of $u$ on $\Omega$. For the $p$-modulus on metric measure spaces
and the properties of upper gradients, see for example
\cite{Hj2}, \cite{HK}, \cite{Sh2}, and \cite{Sh3}.
We use the notation $g_u$ for a $p$-weak upper gradient of $u$. 
The Sobolev space $\sobn\Omega$ consists of those functions $u\in
L^p(\Omega)$ that have a $p$-weak upper gradient $g_u\in L^p(\Omega)$
in $\Omega$. The space $\sobn\Omega$ is a Banach space with the norm
\[
\|u\|_{\sobn\Omega}=\Bigr(\int_\Omega|u|^p\,d\mu
 +\inf_g\int_\Omega|g|^p\,d\mu\Bigl)^{1/p},
\]
where the infimum is taken over all $p$-weak upper gradients $g\in
L^p(\Omega)$ of $u$. 
In the Euclidean space with the Lebesgue measure,
$\sobn{\Omega}=\sob{\Omega}$ for all domains $\Omega\subset\rn$ and
$g_u=|\nabla u|$ is a minimal upper gradient of $u$.

For a measurable set $E\subset X$, the {\it Sobolev space with
zero boundary values} is 
\[
\sobno{E}=\bigl\{u|_E:u\in\sobn{X}\text{ and }u=0
                \text{ in }X\setminus E\bigr\}.
\]
By \cite[Theorem 4.4]{Sh3}, also the space $\sobno E$, equipped with the
norm inherited from $\sobn X$, is a Banach space. 
Note that often the definition of $\sobno{\Omega}$ is given so that
the functions are only required to vanish in 
$X\setminus E$ outside a set of zero $p$-capacity. 
However, our definition gives the same space because functions in
$\sobn X$ are $p$-quasicontinuous by \cite{BBS}.

In order to be able to develop the basic machinery of analysis in the
metric space $X$, we need to assume, in addition to the doubling condition, 
that the geometry of $X$ is rich enough. In practice, this means
that there must exist sufficiently many rectifiable curves everywhere in $X$.
This requirement is in a sense quantified by assuming
that the space $X$ supports a {\it (weak) $(1,p)$-Poincar\'e
inequality}. That is, we assume that there exist constants $c_P>0$ and
$\tau\ge1$ such that for all balls $B\sub X$, all locally integrable
functions $u$ and for all $p$-weak upper gradients $g_u$ of $u$, we have 
\[
\vint{B}|u-u_B|\, d\mu
\le c_Pr\Bigl(\;\vint{\tau B}g_u^p\,d\mu\Bigr)^{1/p},
\]
where 
\[
 u_B=\vint{B}u\,d\mu={\mu(B)}^{-1}\int_Bu\,d\mu
\]
is the integral average of $u$ over $B$.

Standard examples of doubling metric spaces supporting Poincar\'e inequalities
include (weighted) Euclidean spaces, compact Riemannian manifolds, metric graphs,
and Carnot--Carath\'eodory spaces. See for instance \cite{HjK} and
\cite{Hj2}, and the references therein, for more extensive lists of
examples and applications.

\subsection{Capacity and fatness}

Let $\Omega\subset X$ be an open set and let $E\subset \Omega$. The
{\it $p$-capacity} of $E$ with respect to $\Omega$ is
\[
\pcap{E,\Omega}=\inf\int_\Omega g_u^p\,d\mu,
\]
where the infimum is taken over all functions $u\in\sobno{\Omega}$
such that $u\vert_E=1$.
If there are no such functions $u$, then $\pcap{E,\Omega}=\infty$. 
Since the norm of an upper gradient does not increase under
truncation, we may assume that $0\le u\le 1$.  
Note also that because functions in $\sobn X$ are $p$-quasicontinuous
by \cite{BBS}, our definition of $p$-capacity agrees with the
classical definition where admissible functions are required to
satisfy $u=1$ in a neighborhood of $E$.

There exists a constant $C>0$ such that the following comparison between 
the $p$-capacity and measure holds for each $1\leq p<\infty$: 
For all balls $B=B(x,r)$ with $0<r<(1/6)\diam X$ and for each $E\subset B$
\begin{equation}\label{cap ie}
\frac{\mu(E)}{Cr^p}\le\pcap{E,2B}\le\frac{C\mu(B)}{r^p}.
\end{equation}
The lower bound can be obtained by considering $(1,p)$-Poincar\'e
inequality for all admissible functions $0\le u\le 1$ for the capacity
$\pcap{E,2B}$ in the ball $3B$.
For more details, see for example  \cite[Lemma 3.3]{B}.

We say that a set $E\subset X$ is {\it (uniformly) $p$-fat}, $1\le
p<\infty$, if there exists a constant $c_0>0$ such that 
\begin{equation}\label{eq: unif fat}
\pcap{E\cap \overline B(x,r),B(x,2r)}
\ge c_0\pcap{\overline B(x,r),B(x,2r)}
\end{equation}
for all $x\in E$ and all $0<r<(1/6)\diam X$. 
Notice that by the double inequality \eqref{cap ie}, 
$\pcap{\overline B(x,r),B(x,2r)}$ is always comparable to $\mu(B)r^{-p}$.
There are many natural examples of uniformly $p$-fat sets.
For instance, all nonempty subsets of $X$ are uniformly $p$-fat for all $p>s$, 
where $s$ is the doubling dimension of $X$. Also complements of simply
connected subdomains of $\re^2$ and sets satisfying measure density condition 
\[
\mu(B(x,r)\cap E)\ge C\mu(B(x,r))\quad\textrm{for all }x\in E,\quad r>0,
\]
are uniformly $p$-fat for all $1\leq p<\infty$. 
The $p$-fatness condition is stronger than the Wiener
criterion and it is important for example in the study of
boundary regularity of $\mathcal A$-harmonic functions, see \cite{HKM}.

As mentioned in the introduction, uniform fatness is closely related to 
pointwise Hardy inequalities. 

\begin{definition}\label{def: pointwise hardy}
Let $1\le p<\infty$. An open set $\Omega\subsetneq X$ admits {\it the
pointwise $p$-Hardy inequality} if there exist constants $c_H>0$ and
$L\ge1$ such that, for all $u\in\sobno{\Omega}$,
\begin{equation}\label{pointwise m}
|u(x)|\le c_H\dic{x}\Bigl(\M_{L\dic{x}}g_u^p(x)\Bigr)^{1/p}  
\end{equation}
holds at almost every $x\in\Omega$. 
\end{definition}
Above 
\[
\M_R u(x)=\sup_{0<r\le R}\,\vint{B(x,r)}|u|\,d\mu
\]
is the restricted Hardy--Littlewood maximal function of a
locally integrable function $u$. By the maximal theorem \cite[Thm 14.13]{HjK},
$\M_R$ is bounded on $L^p(X)$ for each $1<p\leq\infty$.
Contrary to the Euclidean case, here $\dic x=\dist(x,\Omega^c)$ is the
distance from $x\in\Omega$ to the complement $\Omega^c=X\setminus\Omega$.
We use the this distance because in metric spaces $\dist(x,\bd\Omega)$
may be larger than $\dist(x,\Omega^c)$. 
See however the end of Section \ref{sect: contents} for a related discussion.

\subsection{Main results}

We are now ready to give the general formulation of our main result,
which shows, even in the metric setting, the equivalence between uniform
$p$-fatness of the complement, validity of the pointwise $p$-Hardy
inequality, and two Poincar\'e type inequalities. Here $\tau\geq 1$ is
the dilatation constant from the $(1,p)$-Poincar\'e inequality.
\begin{theorem}\label{thm: eqv}
Let $1\le p<\infty$ and let 
$X$ be a complete, doubling metric measure
space supporting a $(1,p)$-Poincar\'e inequality. Then,
for an open set $\Omega\subsetneq X$,
the following assertions are quantitatively equivalent:
\begin{enumerate}
\item[(a)] The complement $\Omega^c$ is uniformly $p$-fat.
\item[(b)] For all $B=B(w,r)$, with $w\in \Omega^c$
and $r>0$, and every $u\in \sobno{\Omega}$
\begin{equation}\label{eq: bdry poinc}
\int_B |u|^p\,d\mu\le C r^p \int_{5\tau B} g_u^p\,d\mu\,. 
\end{equation}
\item[(c)] For all $x\in\Omega$ and every $u\in \sobno{\Omega}$
\begin{equation}\label{eq: avg pw}
|u_{B_x}|^p\le C \dic{x}^p \vint{20\tau B_x}g_u^p\,d\mu,
\end{equation}
where $B_x=B(x,\dic{x})$.
\item[(d)]The open set $\Omega$
admits the pointwise $p$-Hardy inequality \eqref{pointwise m},
and we may choose the dilatation constant to be $L=20\tau$.
\end{enumerate}
\end{theorem}

\begin{remark}
It can be seen from the proof of Theorem \ref{thm: eqv} that the
conditions (a)--(d) are equivalent also in a local sense, if
interpreted correctly. Indeed, if one of the conditions holds
near a point $w_0\in \Omega^c$, then the other conditions hold
near $w_0$ as well if we only consider sufficiently small radii in the
uniform fatness condition \eqref{eq: unif fat} and
in the Poincar\'e type inequality \eqref{eq: bdry poinc}.
\end{remark}

As the self-improvement of uniform fatness was generalized to the metric
space setting by Bj\"orn, MacManus and
Shanmugalingam in~\cite{BMS}, we obtain for $1<p<\infty$ the following
important corollary to Theorem~\ref{thm: eqv}.

\begin{corollary}\label{coro: self-impro}
For $1<p<\infty$ each of the assertions in Theorem \ref{thm: eqv}
possesses a self-improvement property. 
More precisely, if one of the assertions (a)--(d) holds
for $1<p<\infty$, then there exists some $1<q<p$ so that
the same assertion (and thus each of them) holds with the exponent $q$
and constants depending only on $p$ and the associated data. 
\end{corollary}

Notice that we only assume that $X$ supports a $(1,p)$-Poincar\'e
inequality, but in the above corollary we actually need that $X$
supports a $(1,q)$-Poincar\'e inequality for some $q<p$ as well.
By a result of Keith and Zhong \cite{KeZ}, this is in fact always true if $X$ 
is complete, doubling and supports a weak $(1,p)$-Poincar\'e inequality.

In the previous literature concerning pointwise Hardy inequalities
(see e.g.\ \cite{Hj3} and \cite{Le1}), a sort of a self-improvement has
actually been an a priori assumption when the passage from
pointwise inequalities to the usual Hardy inequality was considered. 
Now, by Corollary \ref{coro: self-impro}, such an extra assumption
becomes unnecessary. 
Especially, using the maximal theorem for an exponent $1<q<p$,
for which $\Omega$ still admits the pointwise inequality, we obtain the
following corollary just as in the Euclidean case.

\begin{corollary}\label{coro: pw implies int}
If an open set $\Omega\sub X$ admits the pointwise $p$-Hardy
inequality \eqref{pointwise m} for some $1<p<\infty$,
then $\Omega$ admits the $p$-Hardy inequality, that is,
there exists $C>0$ such that
\begin{equation*}
\int_\Omega \frac{u(x)^p}{\dic x^p}\,d\mu 
\leq C\int_\Omega g_u(x)^p\,d\mu
\end{equation*}
for every $u\in \sobno{\Omega}$.
\end{corollary}

However, the result of Corollary~\ref{coro: pw implies int}, when
viewed as a consequence of Theorem~\ref{thm: eqv}, depends on a heavy
machinery of non-trivial results already in the Euclidean setting, let
alone in general metric spaces, as the self-improvement of uniform
fatness is involved. In particular, the theory of Cheeger
derivatives is needed in the metric case.
The ideas of Wannebo~\cite{W} lead to an alternative proof for
Corollary~\ref{coro: pw implies int}, which is based on completely 
elementary tools and methods, and especially avoids the use of the
self-improvement. Using this approach, we give in 
Theorem~\ref{thm: wannebo} a direct proof for the fact that 
uniform $p$-fatness of the complement of $\Omega$ implies that $\Omega$ admits 
 the $p$-Hardy inequality. Note that this result was first generalized to 
metric spaces in~\cite{BMS}, but there the proof was based on the self-improvement. 
As the pointwise $p$-Hardy inequality implies the uniform fatness of
the complement by Theorem~\ref{thm: eqv}, Corollary~\ref{coro: pw implies int} 
follows. 

It would also be interesting to acquire an alternative proof for 
Corollary~\ref{coro: self-impro} by showing the self-improvement
directly for one of the conditions (b)--(d) in Theorem~\ref{thm: eqv}.
Let us remark here that self-improving properties of integral Hardy
inequalities were considered in \cite{KZ}, but these results and
methods do not seem apply for pointwise inequalities.

\section{From fatness to pointwise Hardy}\label{sect: from fat to pw}

This section deals with the proofs of the implications
(a)$\Rightarrow$(b)$\Rightarrow$(c)$\Rightarrow$(d) of Theorem \ref{thm: eqv}.
The implication (a)$\Rightarrow$(d), that uniform $p$-fatness of
the complement implies the pointwise $p$-Hardy inequality,
is a generalization of an Euclidean result of Kinnunen and Martio 
\cite[Thm 3.9]{KiMa2} and Haj\l asz \cite[Thm 2]{Hj3}. 

Our proof utilizes the following Sobolev type inequality, proved in
the classical case by Maz'ya (c.f. \cite[Ch.\ 10]{M}) and in the
metric setting by Bj\"orn \cite[Proposition 3.2]{B}. 
We recall the main ideas of the proof for the sake of completeness.
\begin{lemma}\label{lemma: sobo mazya}
There is a constant $C>0$ such that for each $u\in\sobn{X}$ and for
all balls $B\sub X$ we have 
\begin{equation}\label{eq: mazya}
\vint{B}|u|^p\,d\mu\le\frac C{\pcap{\frac 1 2 B\cap\{u=0\},B}}
\int_{5\tau B}g_u^p\,d\mu,
\end{equation}
where $\tau$ is from the $(1,p)$-Poincar\'e inequality.
\end{lemma}

\begin{proof}
Let $B=B(x,r)$ be a ball and let $\vphi$ be a $2/r$-Lipschitz function
such that $0\le\vphi\le1$, $\vphi=1$ on $1/2B$ and $\vphi=0$ outside
$B$. We may assume that $u\ge0$ in $B$. 
The function 
\[
v=\vphi(1-u/\bar u),
\]
where $\bar u=(\,\vint Bu^p\,d\mu)^{1/p}$, is a test function for the
capacity in \eqref{eq: mazya}. The claim follows by estimating the
integral of $g_v^p$, 
\[
g_v=|1-u/\bar u|2r^{-1}+g_u/\bar u.
\]
Here one needs  a $(p,p)$-Poincar\'e inequality, which by
\cite[Theorem 5.1]{HjK} follows from the $(1,p)$-Poincar\'e inequality
with dilatation constant $5\tau$.
\end{proof}

We also need the following pointwise inequality for $N^{1,p}$-functions in
terms of the maximal function of the $p$-weak upper gradient: 
There is a constant $C>0$, depending only on the doubling
constant and the constants of the Poincar\'e inequality, such that 
\begin{equation}\label{umg}
|u(x)-u_{B}|\le Cr\bigl(\M_{\tau r}g_u^p(x)\bigr)^{1/p}
\end{equation}
whenever $B=B(x,r)$ is a ball and $x$ is a Lebesgue point of $u$.
Estimate \eqref{umg} follows easily from 
a standard telescoping argument, see for example \cite{HjKi}.
Note that $u$ has Lebesgue points almost everywhere
in the $p$-capacity sense, see \cite{KL}, \cite{KKST}.

\begin{proof}[Proof of Theorem \ref{thm: eqv} 
(a)$\implies$(b)$\implies$(c)$\implies$(d)]
$ $

\noindent{\bf (a)$\implies$(b):}  
Let $u\in\sobno\Omega$ and let $B=B(w,r)$, 
where $w\in \Omega^c$. Assume first that $0<r < (1/6)\diam X$. 
Since $u$ vanishes outside $\Omega$, we have $\Omega^c\subset\{u=0\}$.
Using the $p$-fatness of $\Omega^c$,
estimate \eqref{cap ie}, and the doubling property of $\mu$, we obtain
\begin{align*} 
\pcap{\tfrac12B\cap\{u=0\},B}
&\ge \pcap{\tfrac12B\cap\Omega^c,B}\\
&\ge c_0\pcap{\tfrac12B,B}
\ge C\mu(B)r^{-p}.
\end{align*} 
This, together with Lemma \ref{lemma: sobo mazya}, gives 
\[
\int_B|u|^p\,d\mu 
\le \frac {C\mu(B)}{\pcap{\tfrac12B\cap\{u=0\},B}}\int_{5\tau B}g_u^p\,d\mu
\le Cr^p\int_{5\tau B}g_u^p\,d\mu.
\]
If $(1/6)\diam X \leq r \leq \diam X$, we take
$\widetilde B=B\big(w,(1/7)\diam X\big)$.
From the triangle inequality it follows that
\[
\int_B|u|^p\,d\mu \leq C\Big(\int_B|u-u_B|^p\,d\mu
+  \mu(B)|u_{\widetilde B}|^p + \mu(B)|u_{\widetilde B}-u_B|^p \Big).
\]
We can then use the $(1,p)$-Poincar\'e inequality, the above case for
the ball $\widetilde B$, and the doubling property,
and the claim for $B$ follows with simple calculations.

Finally, if $r>\diam X$, the claim is clear by the previous cases.

\noindent{\bf (b)$\implies$(c):}  
Let $u\in\sobno\Omega$, $x\in\Omega$, and let $B_x=B(x,\dic x)$. 
Choose a point $w\in \Omega^c$ so that 
\[
R=\dist(x,w)\le2\dic{x},
\]
and let $B_0=B(w,R)$.
Now 
\[
|u_{B_x}|\le|u_{B_x}-u_{B_0}|+|u|_{B_0},
\]
where, by the $(1,p)$-Poincar\'e inequality, the fact that
$B_0\subset 4B_x$ and $B_x\subset2B_0$, and the doubling property,
\[
|u_{B_x}-u_{B_0}|\le C \dic x\Bigl(\;\vint{4\tau B_x}g_u^p\,d\mu\Bigr)^{1/p}.
\] 
Using the H\"older inequality, assumption (b), and the doubling
property, we obtain
\[
|u|_{B_0}
\le\Bigl(\;\vint{B_0}|u|^p\,d\mu\Bigr)^{1/p}
\le CR\Bigl(\;\vint{5\tau B_0}g_u^p\,d\mu\Bigr)^{1/p}
\le C\dic x\Bigl(\;\vint{20\tau B_x}g_u^p\,d\mu\Bigr)^{1/p}.
\]
The claim follows by combining these two estimates.

\noindent{\bf (c)$\implies$(d):}  
Let $u\in\sobno\Omega$ and let $x\in\Omega$ be a Lebesgue point of $u$.  
Now 
\[
|u(x)|\le|u(x)-u_{B_x}|+|u_{B_x}|,
\]
where, by \eqref{umg} 
\[
|u(x)-u_{B_x}|\le C\dic x\bigl(\M_{\tau\dic x}g_u^p(x)\bigr)^{1/p},
\]  
and by (c)
\[
|u_{B_x}|\le C\dic x\Bigl(\;\vint{20\tau B_x}g_u^p\,d\mu\Bigr)^{1/p}
\le C\dic x\bigl(\M_{20\tau\dic x}g_u^p(x)\bigr)^{1/p}.
\]
The pointwise $p$-Hardy inequality follows from the above estimates.
\end{proof}

By slightly modifying the proof above or the proof in 
\cite[Thm 3.9]{KiMa2}, we obtain a $p$-Hardy inequality containing a
fractional maximal function of the upper gradient.
\begin{corollary}\label{cor: phardy from pfat}
Let $1\le p< \infty$ and let $\Omega\sub X$ be an open set
whose complement is uniformly $p$-fat. 
Then there is a constant $C>0$, independent of $\Omega$, such that
for all $0\le\alpha<p$ and for all $u\in\sobno{\Omega}$,
\begin{equation}\label{phardy from pfat alpha}
|u(x)|\le C\dic{x}^{1-\alpha/p}\bigl(\M_{\alpha,20\tau\dic{x}}g_u^p(x)\bigr)^{1/p}
\end{equation} 
whenever $x\in\Omega$ is a Lebesgue point of $u$.
\end{corollary}
Here, for $\alpha\ge0$, the restricted fractional maximal function of
a locally integrable function $u$ is 
\[
\M_{\alpha,R} u(x)=\sup_{0<r\le R}\,r^{\alpha}\vint{B(x,r)}|u|\,d\mu.
\]

\section{From pointwise Hardy to fatness}\label{sect: from pw to fat}

In this section we prove the following lemma, from which 
the part (d)$\Rightarrow$(a) of Theorem \ref{thm: eqv}
and the previously unknown necessity part of Theorem \ref{thm: main}
follow.

\begin{lemma}\label{lemma: p-fat from pw}
Let $1\le p<\infty$ and let $\Omega\sub X$ be an open set. 
If $\Omega$ admits the pointwise $p$-Hardy inequality
\eqref{pointwise m},
then $\Omega^c$ is uniformly $p$-fat. The constant in the uniform
fatness condition \eqref{eq: unif fat}
depends only on $p$, $c_H$, and the constants related to $X$.
\end{lemma}

\begin{proof}
Let $B=B(w,R)$, where $w\in\Omega^c$ and $0<R<(1/6)\diam X$.
By \eqref{cap ie}, it suffices to find a constant $C>0$,
independent of $w$ and $R$, such that 
\begin{equation}\label{eq: goal}
\mu(B)R^{-p}\leq C\int_{2B} g_v^p\,d\mu
\end{equation}
whenever $g_v$ is an upper gradient of a function $v\in\sobno{2B}$
satisfying $0\leq v\leq 1$ and $v=1$ in $\Omega^c\cap\ol B$. By the
quasicontinuity of $N^{1,p}$-functions, we may assume that $v=1$ in an
open neighborhood of $\Omega^c\cap \ol B$.

Let $l=\big[2(L+1)\big]^{-1}$, where $L$ is from the pointwise
$p$-Hardy inequality \eqref{pointwise m}.
The doubling condition implies that $\mu(lB)\geq l^s \mu(B)/c_D$.
If now $v_B > l^s/2c_D$, we obtain from the Poincar\'e inequality 
for $v\in\sobno{2B}$ (see for example \cite[Proposition 3.1]{B}) that
\[
1 \leq C \vint{B} |v|\,d\mu
\leq C R \bigg(\,\vint{2B} g_v^p\,d\mu\bigg)^{1/p},
\]
and \eqref{eq: goal} follows by the doubling condition. 

We may hence assume that $v_B\leq l^s/2c_D$.
Let $\psi\in\sobno{B}$ be a cut-off function, defined as 
\[
\psi(x) 
= \max\Bigl\{0,1-\tfrac 4 R \dist\bigl(x,\tfrac 1 2 B\bigr)\Bigr\},
\]
and take 
\[u=\min\{\psi,1-v\}.\]
Since $1-v=0$ in an open set
containing $\Omega^c\cap B$ and  
$N^{1,p}(X)$ is a lattice, 
we have that $u\in\sobno{\Omega}$.
Moreover, $u$ has an upper gradient $g_u$ such
that $g_u=g_v$ in $1/2 B$. 

We define $C_1 = l^s/4c_D$ and 
\[
E=\big\{  x\in l B : u(x) > C_1 
\text{ and }\eqref{pointwise m}\text{ holds for }u\text{ at }x\big \},
\]
and claim that
\begin{equation}\label{eq: eka}
 \mu(E) \geq C_1\mu(B).
\end{equation}
To see this, first notice that $u=1-v$ in $lB$ and that  
$\mu(lB)\geq 4C_1 \mu(B)$. As $v_B\leq l^s/2c_D = 2C_1$, we obtain 
\begin{equation}\label{eq: big int}
\begin{split}
\int_{ l B } u\,d\mu  
=\int_{ l B } (1-v)\,d\mu  
& \geq \int_B (1-v)\,d\mu -\mu(B\setminus l B)\\
& \geq (1-2C_1)\mu(B) - \mu(B) +  \mu(lB)\\
& \geq 2 C_1\mu(B).
\end{split}
\end{equation} 
Since the pointwise $p$-Hardy holds for almost every $x\in\Omega$, 
we have $u\leq C_1$ almost everywhere in $lB\setminus E$.
Thus a direct computation using estimate \eqref{eq: big int} yields
\eqref{eq: eka}:
\[
\begin{split}
  \mu(E) 
& \geq \int_E u\,d\mu = \int_{l B} u\,d\mu 
                             - \int_{l B\setminus E} u\,d\mu\\
& \geq 2C_1\mu(B) - \int_{l B} C_1\,d\mu \\
& \geq 2C_1\mu(B) - C_1\mu(B) = C_1\mu(B).
\end{split}
\]

To continue the proof, we fix for each $x\in E$ a radius
$0<r_x\leq L\dic{x}$ such that
\[
 \M_{L\dic{x}}g_u^p(x)\leq 2 \,\vint{B(x,r_x)} g_u^p\,d\mu.
\]
By the standard $5r$-covering theorem (see e.g.\ \cite{CW}), there are
pairwise disjoint balls $B_i=B(x_i,r_i)$, where $x_i\in E$ and
$r_i=r_{x_i}$ are as above, so that $E\sub\bigcup_{i=1}^{\infty} 5
B_i$. It follows immediately from \eqref{eq: eka} and the doubling
condition that
\begin{equation}\label{eq: start}
\mu(B)\leq C_1^{-1} \mu(E) 
\leq C\sum_{i=1}^{\infty}\mu(B_i).
\end{equation}
As $x_i\in lB$ and $w\notin\Omega$, we have $\dic{x_i}\leq lR$. 
Hence, by the choice of $l$, we obtain for each $y\in B_i$ that
\[
\dist(w,y)\leq \dist(w,x_i) + \dist(x_i,y) 
\leq lR + L\dic{x_i}\leq lR(1 + L)= R/2,
\]
and so $B_i\sub1/2 B$. 
This means, in particular, that $g_u=g_v$ in each $B_i$. 
Since $u(x_i)>C_1$ for each $i$, the pointwise $p$-Hardy inequality 
\eqref{pointwise m} and the choice of the radii $r_i$ imply that
\[
C_1^p 
\leq |u(x_i)|^p\leq C \dic{x_i}^p \M_{L\dic{x_i}}g_u^p(x) 
\leq C R^p \mu(B_i)^{-1} \int_{B_i}g_u^p\,d\mu,
\]
and so
\[
\mu(B_i)\leq C R^p \int_{B_i}g_v^p\,d\mu.
\]
Inserting this into \eqref{eq: start} leads us to
\[
\mu(B)\leq C R^p \sum_{i=1}^{\infty} \int_{B_i} g_v^p\, d\mu
\leq C R^p \int_{2B} g_v^p\,d\mu,
\]
where we used the fact that the balls $B_i\sub 2B$ are pairwise disjoint.
This proves estimate \eqref{eq: goal}, and the lemma follows.
\end{proof}

\section{From fatness to Hardy}\label{sect: from fat to int}

The purpose of this section is to give a straight-forward proof for the
fact that uniform $p$-fatness of the complement $\Omega^c$
suffices for $\Omega$ to admit the $p$-Hardy inequality. 
Our proof follows the ideas of Wannebo~\cite{W}. 
A similar method was also used in \cite{BK} in the context of
Orlicz--Hardy inequalities. As mentioned earlier, the following
result first appeared in the metric space setting in~\cite{BMS},
where the proof was based on the self-improvement of uniform $p$-fatness.

\begin{theorem}\label{thm: wannebo}
Let $1<p<\infty$ and let $\Omega\subset X$ be an open set. 
If $\Omega^c$ is uniformly $p$-fat then $\Omega$ admits the $p$-Hardy
inequality, quantitatively.
\end{theorem}

\begin{proof}
To make the proof as simple as possible, let us assume that the dilatation
constant in the right-hand side of Theorem \ref{thm: eqv} (b) is 2. The
general case follows by obvious modifications. 
Let
\[
 \Omega_n=\{x\in\Omega\,:\, 2^{-n}\leq\dic x<2^{-n+1}\} 
\]
and 
\[
 \widetilde\Omega_n=\bigcup_{k=n}^\infty\Omega_k.
\]
Let $\Fc_n$ be a cover of $\Omega_n$ with balls of radius
$2^{-n-2}$ such that their center points are not included in any other
ball in $\Fc_n$. Associate to each ball $B\in \Fc_n$ a bigger ball 
$\widetilde B\supset B$, whose radius is $2^{-n+2}$ and whose
center point is on $\bd\Omega$. 
Note that $2\widetilde B\cap \Omega \subset \widetilde\Omega_{n-2}$ and that
\[
\sum_{B\in\Fc_n}\chi_{B}<C\quad\textrm{ and }\quad 
\sum_{ B\in {\Fc_n}}\chi_{2\widetilde B}<C,
\]
where the constant $C>0$ only depends on the doubling constant of $\mu$.

Let $u\in\sobno\Omega$. The condition (b) of Theorem \ref{thm: eqv}
(which follows from the uniform $p$-fatness of the complement) 
implies that for every $B\in\Fc_n$ we have
\[
 \int_{B}|u|^p\,d\mu\leq\int_{\widetilde B}|u|^p\,d\mu
\leq C 2^{-np}\int_{2\widetilde B}g_u^p\,d\mu.
\]
By summing up the inequalities above, we obtain
\begin{equation}\label{eqn:summattavat}
\begin{split}
\int_{\Omega_n}|u|^p d\mu
&\leq \sum_{B\in\Fc_n}\int_{B}|u|^pd\mu
\leq C 2^{-np}\sum_{ B\in{\Fc_n}} \int_{2\widetilde B}g_u^pd\mu\\
&\leq C 2^{-np}\int_{\widetilde\Omega_{n-2}}g_u^pd\mu.
\end{split}
\end{equation}
Let $0<\beta<1$ be a small constant to be fixed later. We
multiply~\eqref{eqn:summattavat} by $2^{n(p+\beta)}$ and sum the
inequalities to obtain
\begin{equation}\label{eqn:summaus}
\begin{aligned}
\int_\Omega |u(x)|^p\dic x^{-p-\beta}\, d\mu
&\leq \sum_{n=-\infty}^\infty \int_{\Omega_n}|u(x)|^p2^{n(p+\beta)}\,d\mu\\
&\leq C\sum_{n=-\infty}^\infty2^{n\beta}\int_{\widetilde\Omega_{n-2}}g_u(x)^p\,d\mu \\
&=C\sum_{k=-\infty}^\infty \biggl(\,\sum_{n=0}^\infty 2^{(k+2-n)\beta}
  \int_{\Omega_k}g_u(x)^p\,d\mu\biggr)\\
& \leq C \sum_{k=-\infty}^\infty \frac{2^{k\beta}}{\beta}\int_{\Omega_k}g_u(x)^p\,d\mu\\
&\leq\frac{C}{\beta}\int_\Omega g_u(x)^p\operatorname{d}_\Omega(x)^{-\beta}\,d\mu.
\end{aligned}
\end{equation}
In the calculations above, we used the fact that 
$2^{-k}\leq \dic x\leq 2\cdot 2^{-k}$ for every $x\in\Omega_k$.

Now let $v\in\sobn{\Omega}$ be a function with a compact support in
$\Omega$ and let
\[
u(x)=v(x)\dic x^{\beta/p}.
\]
Then the function 
\[
g_u(x)= g_v(x) \dic x^{\beta/p}+ \frac{\beta}{p} v(x) \dic x^{\beta/p-1}
\]
is a $p$-weak upper gradient of $u$. 
Thus, by~\eqref{eqn:summaus}, we have
\[
\begin{split}
 \int_\Omega \frac{v(x)^p}{\dic x^p}\,d\mu 
&= \int_\Omega \frac{u(x)^p}{\dic x^{p+\beta}}\,d\mu
\leq \frac{C}{\beta}\int_\Omega \frac{g_u(x)^p}{\dic x^\beta}\,d\mu\\
&\leq \frac{C}{\beta}\int_\Omega g_v(x)^p\,d\mu
+\frac{C}{\beta}\frac{\beta^p}{p^p}\int_\Omega \frac{v(x)^p}{\dic x^p}\,d\mu.
\end{split}
\]
If $\beta>0$ is small enough, the last term on the right-hand side
can be included on the left-hand side and we obtain
\[
\int_\Omega \frac{v(x)^p}{\dic x^p}\,d\mu 
\leq C\int_\Omega g_v(x)^p\,d\mu.
\]
This completes the proof because functions with compact support are
dense in $\sobno\Omega$, see \cite[Theorem 4.8]{Sh3}.
\end{proof}

Notice that the requirement $p>1$ is essential in Theorem \ref{thm: wannebo}.
For instance, smooth domains in $\R^n$ admit the pointwise $1$-Hardy inequality but not the integral $1$-Hardy.

\section{Hausdorff contents}\label{sect: contents}
It is well-known that capacities and Hausdorff contents are closely related both
in Euclidean spaces and general metric spaces, see e.g.\ \cite{HKM} and \cite{HK}.
In metric spaces we follow
\cite{A}, \cite{AMP}, and \cite{KKST}, and use Hausdorff contents $\Hh^t_R$,
defined by applying the Carath\'eo\-do\-ry construction to functions 
\[
h(B(x,r))=\frac{\mu(B(x,r))}{r^t},
\]q
where $r\leq R$.
Thus the {\it Hausdorff content of codimension $t$} of a set $E\subset
X$ is given by  
\[
\Hh^t_R(E)=\inf\Bigl\{\sum_{i\in I} h(B(x_i,r_i)):
E\subset\bigcup_{i\in I} B(x_i,r_i),\ r_i\leq R \Bigr\}.
\]
Here we may actually assume that $x_i\in E$, as this increases $\Hh^t_R(E)$
at most by a constant factor.
 
If the space $X$ is $Q$-regular, then $\Hh^t_\infty(E)$ is comparable
with the usual Hausdorff content $\Hc^{Q-t}_{\infty}(E)$, which is
defined by using the gauge function $h(B(x,r))={r^{Q-t}}$. Recall that
$Q$-regularity means that there are constants $c_1,c_2>0$ such that
\[
c_1r^Q\le\mu(B(x,r))\le c_2r^Q
\]
for all balls $B(x,r)$ in $X$.

Now, by slightly modifying the argument in \cite[Thm.\ 5.9]{HK} 
(see also \cite{Co} and \cite{DGP}),
one can show that if $E\sub X$ is a closed set and 
there exists some $1\leq q<p$ and a constant $C>0$ so that
for all $w\in E$ and every $R>0$,
\begin{equation}\label{eq: c to a}
\Hh_{R/2}^q\big(E \cap \ol B(w,R)\big)\geq C \mu(\ol B(w,R))R^{-q}, 
\end{equation}
then $E\sub X$ is uniformly $p$-fat. 
Conversely, by rewriting the argument of 
\cite[Thm.\ 2.27]{HKM} (see also \cite[Thm.\ 4.9]{Co})
for the content $\Hh_{R/2}^p$, it is not hard to see that uniform
$p$-fatness of $E$ leads to \eqref{eq: c to a}, but with $q$ replaced
by $p$. Using the self-improvement of uniform fatness, we then conclude that
uniform $p$-fatness of $E$ implies the existence of an exponent $q<p$ 
for which \eqref{eq: c to a} holds.
Hence \eqref{eq: c to a}, with an exponent $1\leq q <p$,
is actually equivalent with the uniform $p$-fatness of $E$. 

In this section, we investigate similar density conditions for the
boundary of a domain $\Omega\sub X$. To this end, we consider a
version of the pointwise Hardy inequality where the distance is taken
to the boundary instead of the complement. We define
\[
\dom(x)=\dist(x,\bd\Omega)\text{ for }x\in\Omega. 
\]
The following lemma is a metric space 
generalization of a result from \cite{Le1} and \cite{DGP}.

\begin{lemma}\label{thm: inner density}
Let $1\le p< \infty$ and let $\Omega\subset X$ be an open set. 
Assume that $\Omega$ admits the pointwise $p$-Hardy inequality
\begin{equation}\label{bdry pointwise m}
|u(x)|\le c_H\dom(x)\Bigl(\M_{L\dom(x)}g_u^p(x)\Bigr)^{1/p}
\end{equation}
for all $u\in\sobno{\Omega}$.
Then 
\begin{equation}\label{eq: inner density}
\Hh^{p}_{\dom(x)}\big(\bd\Omega\cap \ol B(x,2L\dom(x))\big)
\geq C \dom(x)^{-p}\mu\big(\ol B(x,\dom(x))\big).
\end{equation}
for all $x\in\Omega$.
\end{lemma}

\begin{proof}
Let $x\in\Omega$. We define $R=\dom(x)$, $B=\ol B(x,R)$, and
\[
E=\bd\Omega\cap 2LB.
\]
Let $\{B_i\}_{i=1}^N$, where $B_i=B(w_i,r_i)$ with $w_i\in E$ and $0<r_i\leq R$, 
be a covering of $E$; we may assume that the covering is finite by the
compactness of $E$.

It is now enough to show that there exists a constant $C>0$,
independent of the particular covering, such that
\begin{equation}\label{eq: sum is big}
\sum_{i=1}^N \mu(B_i){r_i}^{-p} \geq C \mu(B) {R}^{-p}. 
\end{equation}
If $r_i\geq R/4$ for some $1\leq i\leq N$, then, by
\eqref{doubling dimension} and the fact that $r_i^{-p}\geq R^{-p}$, we have
\[
\mu(B_i)r_i^{-p}\ge C 
\mu(B)\Bigl(\frac {r_i}R\Bigr)^{s}R^{-p}\ge C 
\mu(B) R^{-p},
\]
from which \eqref{eq: sum is big} readily follows. 

We may hence assume that $r_i < R/4$ for all $1\leq i\leq N$.
Now, define
\[
\vphi(y)=\min_{1\leq i\leq N}\big\{1,\,r_i^{-1}\dist(y,B_i)\big\}\,
\]
and let $\psi\in \sobno{2LB}$ be a cut-off function 
such that $0\leq\psi\leq 1$ and $\psi(y)=1$ for all $y\in LB$.
Then the function
\[
u=\min\{\psi,\vphi\}\ch{\Omega}
\]
belongs to $\sobno{\Omega}$.
As $r_i<R/4$ for all $1\leq i\leq N$, it follows that
$\dist(x,2B_i)\ge R/2$ for all $1\leq i\leq N$, and so $u(x)=1$.

In addition, $u$ has an upper gradient $g_u$ such that
\begin{equation}\label{eq: grad}
{g_u(y)}^p\leq \sum_{i=1}^N r_i^{-p}\ch{2B_i}(y)
\end{equation}
for a.e.\ $y\in LB$.
Especially, we must have $r>R/2$ in order to obtain something positive
when estimating $\M_{LR}{g_u}^p(x)$.
As the pointwise inequality \eqref{bdry pointwise m}
holds for the continuous function $u\in\sobno{\Omega}$
at every $x\in\Omega$, we have 
\[
\begin{split}
1 & =  |u(x)|^p 
\leq CR^p \M_{LR}g_u^p(x) 
\leq C {R}^{p} \sup_{R/2 \leq r \leq LR}\;\vint{B(x,r)} g_u^p\,d\mu\\
    & \leq C  {R}^{p} \mu\big(\tfrac 1 2 B\big)^{-1}
        \int_{LB} g_u^p\,d\mu
\leq  C {R}^{p} \mu(B)^{-1} 
       \sum_{i=1}^N \mu(2B_i)\, r_i^{-p},
\end{split}
\]
where the last inequality is a consequence of \eqref{eq: grad}. 
Estimate \eqref{eq: sum is big} then 
easily follows with the help of the doubling property.
\end{proof}

Next we show that the inner boundary density condition \eqref{eq: inner density}
is actually almost equivalent to the pointwise $p$-Hardy
inequality. The proof below uses an idea from \cite{HK}, but is new in
the context of Hardy inequalities.

\begin{theorem}\label{thm: pw from dense}
Let $1<p<\infty$ and let $\Omega\sub X$ be an open set. If estimate
\eqref{eq: inner density} holds with an exponent $1\leq q<p$ for all
$x\in\Omega$, then $\Omega$ admits the pointwise $p$-Hardy inequality
\eqref{bdry pointwise m}, but possibly with a different dilatation
constant in the maximal function.
\end{theorem}

\begin{proof}
Let us first assume that $u\in \sobno{\Omega}$ has a compact support
in $\Omega$. Let $B=\ol B(x,R)$, where $x\in \Omega$ and $R=\dom(x)$. 
We are going to show that 
\begin{equation}\label{eq: avg pw bdry}
|u_{B}|^p\leq C \dom(x)^p \vint{3\tau LB}g_u^p\,d\mu,
\end{equation}
where $C>0$ and $\lambda\geq 1$ are independent of $x$,
whence the pointwise $p$-Hardy inequality follows for almost every $x\in\Omega$ by
 Theorem \ref{thm: eqv}.

If $u_B=0$, the claim \eqref{eq: avg pw bdry} is true, and so we may
assume that $|u_B|>0$, and in fact, by homogeneity, that $|u_B|=1$.
Let $w\in \bd\Omega\cap 2L B$ and let $B_k=B(w,r_k)$, where 
$r_k=(5 \tau 2^k)^{-1}R$, $k\in\N$.
It then follows that
\[
1 = |u(w)-u_B|\leq |u_{B_0}|+|u_{B_0}-u_B|.
\]
Now, if $|u_{B_0}|<1/2$, we infer, using
the $(1,p)$-Poincar\'e inequality, the facts $B_0\sub 3LB$ and $B\sub 3LB_0$, 
and the doubling property, that
\[
 \tfrac 1 2 
\leq |u_{B_0}-u_B|\leq |u_{B_0}-u_{3B}|+|u_{B}-u_{3B}|
\leq C R \Big(\, \vint{3\tau LB} g_u^p\,d\mu\Big)^{1/p}. 
\]
As $|u_B|=1$, the claim follows.
 
Thus we may assume that $1/2 \leq |u_{B_0}|=|u(w)-u_{B_0}|$
for every $w\in \bd\Omega\cap 2L B$.
A standard chaining argument, using the $(1,p)$-Poincar\'e inequality
(see for example \cite{HjK}) and the assumption that the support of $u$ is compact,
leads us to estimate
\begin{equation}\label{eq: chain}
 1\leq C \sum_{k=0}^\infty r_k\, \Bigl(\;\vint{\tau B_k}g_u^p\,d\mu\Bigr)^{1/p}.
\end{equation}
From \eqref{eq: chain} it follows that there must be a constant
$C_1>0$, independent of $u$ and $w$, and at least one index $k_w\in\N$
such that
\[
r_{k_w}\, \Bigl(\;\vint{\tau B_{k_w}}g_u^p\,d\mu\Bigr)^{1/p} \geq 
   C_1 2^{-k_w(1-q/p)} = C_1 \left(\frac{r_{k_w}}{R}\right)^{1-q/p}.
\]
In particular, we obtain for each $w\in \bd\Omega\cap 2LB$ a radius 
$r_w\le R/(5\tau)$ and a ball $B_w=B(w,r_w)$ such that
\begin{equation}\label{eq: good r}
\mu(\tau B_w)\,r_w^{-q}\leq C R^{p-q}\int_{\tau B_w}g_u^p \,d\mu.
\end{equation}
Again, the $5r$-covering lemma implies the existence of points 
$w_1,w_2,\ldots,w_N\in \bd\Omega\cap 2LB$ such that
if we set  $r_i=r_{w_i}$, then the balls $\tau B_i=B(w_i,\tau r_i)$
are pairwise disjoint, but still
$\bd\Omega\cap 2LB\sub \bigcup_{i=1}^N 5\tau B_i$. 
Assumption \eqref{eq: inner density}, the doubling property, 
estimate \eqref{eq: good r}, and the pairwise disjointness of the
balls $\tau B_i\sub 3\tau LB$ then yield
\begin{equation}\label{eq: Hastimate}
\begin{split}
R^{-q}\mu(B)
& \leq C\Hh^q_{R}(\bd\Omega\cap 2LB)\\
& \leq C\sum_{i=1}^N \mu(5\tau B_{i})(5\tau r_{i})^{-q} 
\leq C \sum_{i=1}^N \mu(\tau B_{i}){r_{i}}^{-q}\\
& \leq C\sum_{i=1}^N R^{p-q}
   \int_{\tau B_{i}}g_u^p \,d\mu
 \leq C R^{p-q}
   \int_{3\tau LB}g_u^p \,d\mu.
\end{split}
\end{equation}
As we assumed $|u_B|=1$, estimate \eqref{eq: avg pw bdry} now follows
from \eqref{eq: Hastimate} and the doubling condition.

For a general $u\in \sobno{\Omega}$ estimate \eqref{eq: avg pw bdry}
follows by a suitable approximation with compactly supported
functions. 
\end{proof}

If there now exists a constant $C\geq 1$ such that
\begin{equation}\label{eq: distances}
\dic{x}\leq\dom(x)\leq C\dic{x}\ \text{ for each } x\in\Omega, 
\end{equation}
then it is clear that pointwise inequalities \eqref{pointwise m}
and \eqref{bdry pointwise m} are quantitatively equivalent. In
particular, if the inner boundary density condition \eqref{eq: inner density} 
with codimension $q$ holds for all $x\in\Omega$, then Theorems
\ref{thm: pw from dense} and \ref{thm: eqv} imply that $\Omega^c$ is
uniformly $p$-fat for all $p>q$.
On the other hand, easy examples show that $\Omega^c$ need not be
uniformly $q$-fat, or equivalently, $\Omega$ need not admit the
pointwise $q$-Hardy inequality, if $q>1$. Hence some information is
inevitably lost once we pass from the pointwise $p$-Hardy inequality
or uniform $p$-fatness (for $1<p<\infty$) to Hausdorff contents; 
in the case $p=1$ there is indeed an equivalence, cf.\ \cite{KKST}.
However, by the self-improvement of the assertions
of Theorem \ref{thm: eqv}, we can still have the following equivalent
characterization in terms of Hausdorff contents (see also \cite{Le1}
and \cite{DGP}). Note that here we need to use again the fact that $X$
supports a $(1,q)$-Poincar\'e inequality for some $q<p$.

\begin{corollary}\label{coro: in fat is eqv}
Assume that $\Omega\sub X$ is such that \eqref{eq: distances} holds.
Then all of the assertions in Theorem \ref{thm: eqv}, with an exponent
$1<p<\infty$, are (quantitatively) equivalent to the following density
condition: There exist some $1<q<p$ and constants $C>0$ and $L\geq 1$
such that
\[
\Hh^{q}_{\dom(x)}\big(\bd\Omega\cap \ol B(x,L\dom(x))\big)
\geq C \dom(x)^{-q}\mu\big(\ol B(x,\dom(x))\big)
\]
for all $x\in\Omega$.
\end{corollary}

It is worth a mention that uniform $p$-fatness of the boundary $\bd\Omega$ is
of course sufficient for the uniform $p$-fatness of
the complement and the pointwise $p$-Hardy inequality, 
but not necessary, as cusp-type domains in $\R^n$, $n\geq 3$, show (cf.\ \cite{Le1}).
Thus it really is essential that we consider above the 
density of the boundary only 
as seen from within the domain, in the sense of \eqref{eq: inner density}.

\subsection*{Acknowledgments} 
The authors wish to thank Juha Kinnunen for numerous helpful
discussions and valuable suggestions during the preparation of this work.

\vspace{0.3cm}
\noindent
\small{\textsc{R.K.},}
\small{\textsc{Department of Mathematics and Statistics},}
\small{\textsc{P.O. Box 68},}
\small{\textsc{FI-00014 University of Helsinki},}
\small{\textsc{Finland}}\\
\footnotesize{\texttt{ riikka.korte@helsinki.fi}}

\vspace{0.3cm}
\noindent
\small{\textsc{J.L.},}
\small{\textsc{Department of Mathematics and Statistics},}
\small{\textsc{P.O. Box 35 (MaD)},}
\small{\textsc{FI-40014 University of Jyv\"askyl\"a},}
\small{\textsc{Finland}}\\
\footnotesize{\texttt{ juha.lehrback@jyu.fi}}

\vspace{0.3cm}
\noindent
\small{\textsc{H.T.},}
\small{\textsc{Department of Mathematics and Statistics},}
\small{\textsc{P.O. Box 35 (MaD)},}
\small{\textsc{FI-40014 University of Jyv\"askyl\"a},}
\small{\textsc{Finland}}\\
\footnotesize{\texttt{ heli.m.tuominen@jyu.fi}}


\begin{thebibliography}{66}
\bibitem{A}{\sc L. Ambrosio},
{\it Fine properties of sets of finite perimeter in doubling metric
  measure spaces}, 
Set-Valued Anal. {\bf 10} (2002), no.2-3, 111--128.
\bibitem{AMP}{\sc L. Ambrosio, M. Miranda, Jr}, and {\sc D. Pallara},
{\it Special functions of bounded variation in doubling metric measure
  spaces},  
In: ``Calculus of variations: Topics from the mathematical heritage of E.\
De Giorgi'', Quad. Mat., 14, Dept. Math., Seconda Univ. Napoli, 2004, 1--45. 
\bibitem{An}{\sc A. Ancona},
{\it On strong barriers and an inequality of Hardy for domains in
  $R\sp n$},
J. London Math. Soc. (2) {\bf 34} (1986), no.2, 274--290.
\bibitem{BBS}{\sc A. Bj\"orn, J. Bj\"orn}, and {\sc N. Shanmugalingam},
{\it Quasicontinuity of Newton--Sobolev functions and density of
  Lipschitz functions in metric measure spaces},
Houston J. Math., {\bf 34} (2008), No.4, 1197--1211. 
\bibitem{B}{\sc J. Bj\"orn},
{\it Boundary continuity for quasiminimizers on metric spaces},
Illinois J. Math. {\bf 46} (2002), no.2, 383--403. 
\bibitem{BMS}{\sc J. Bj\"orn, P. MacManus}, and {\sc N. Shanmugalingam},
{\it Fat sets and pointwise boundary estimates for $p$-harmonic
functions in metric spaces}, 
J. Anal. Math. {\bf 85} (2001), 339--369.
\bibitem{BK}{\sc S.M. Buckley} and {\sc P. Koskela}, 
{\it Orlicz-Hardy inequalities},  
Illinois J. Math. {\bf 48} (2004), no.3, 787--802.
\bibitem{CW}{\sc R.R. Coifman} and {\sc G. Weiss}, 
``Analyse harmonique non-commutative sur certains espaces homog\`enes'',
Lecture Notes in Mathematics, Vol.242. Springer-Verlag, Berlin-New York, 1971. 
\bibitem{Co}{\sc S. Costea},
{\it Sobolev capacity and Hausdorff measures in metric measure spaces},   
Ann. Acad. Sci. Fenn. Math. {\bf 34} (2009), 179--194.
\bibitem{DGP}{\sc D. Danielli, N. Garofalo}, and {\sc N.C. Phuc},
{\it Inequalities of Hardy-Sobolev type in Carnot-Carath\'eodory spaces}, 
In: ``Sobolev Spaces in Mathematics I'', 117-151, Springer, 2009. 
\bibitem{Hj2}{\sc P. Haj\l asz}, 
{\it Sobolev spaces on metric-measure spaces},
In: ``Heat kernels and analysis on manifolds, graphs, and metric spaces'',
(Paris, 2002), 173--218, Contemp. Math. {\bf 338}, Amer. Math. Soc.
Providence, RI, 2003.
\bibitem{Hj3}{\sc P. Haj\l asz}, 
{\it Pointwise Hardy inequalities},  
Proc. Amer. Math. Soc. {\bf 127} (1999), no.2, 417--423.
\bibitem{HjKi}{\sc P. Haj\l asz} and {\sc J. Kinnunen}, 
{\it H\"older quasicontinuity of Sobolev functions on metric spaces},
Rev. Mat. Iberoamericana {\bf 14} (1998), no.3, 601--622.
\bibitem{HjK}{\sc P. Haj\l asz} and {\sc P. Koskela},
{\it Sobolev met Poincar\'e},
Mem. Amer. Math. Soc. {\bf 145} (2000), no.688.
\bibitem{HLP}{\sc G.H. Hardy, J.E. Littlewood}, and {\sc G. P\'olya},
``Inequalities'', (Second edition), 
Cambridge University Press, 1952, Cambridge.
\bibitem{HKM}{\sc J. Heinonen, T. Kilpel\"ainen}, and {\sc O. Martio},
``Nonlinear Potential Theory of Degenerate Elliptic Equations'',
Oxford University Press, 1993, Oxford-New York-Tokyo.
\bibitem{HK}{\sc J. Heinonen} and {\sc P. Koskela}, 
{\it Quasiconformal maps on metric spaces with controlled geometry},
Acta Math. {\bf 181} (1998), 1--61.
\bibitem{KeZ}{\sc S. Keith} and {\sc X. Zhong},
{\it The Poincar\'e inequality is an open ended condition},  
Ann. of Math. (2) {\bf 167} (2008), no.2, 575--599.
\bibitem{KiKiMa}{\sc T. Kilpel\"ainen, \sc J. Kinnunen}, and {\sc O. Martio},
{\it Sobolev spaces with zero boundary values on metric spaces}, 
Potential Anal., {\bf 12} (2000),  no.3, 233--247.
\bibitem{KL}{\sc J. Kinnunen} and {\sc V. Latvala},
{\it Lebesgue points for Sobolev functions on metric spaces},
Rev. Mat. Iberoamericana, {\bf 18} (2002), 685--700.
\bibitem{KiMa2}{\sc J. Kinnunen} and {\sc O. Martio},
{\it Hardy's inequalities for Sobolev functions},
Math. Res. Lett. {\bf 4} (1997), no.4, 489--500.
\bibitem{KKST}{\sc J. Kinnunen, R. Korte, N. Shanmugalingam}, and {\sc H. Tuominen},
{\it Lebesgue points and capacities via boxing inequality in metric
  spaces},  
Indiana Univ. Math. J. {\bf 57} (2008), no.1, 401--430. 
\bibitem{KoSh}{\sc R. Korte} and {\sc N. Shanmugalingam},
{\it Equivalence and self-improvement of $p$-fatness and Hardy's
  inequality, and association with uniform perfectness},
to appear in Math. Z.
\bibitem{KZ}{\sc P. Koskela} and {\sc X. Zhong},
{\it Hardy's inequality and the boundary size},  
Proc. Amer. Math. Soc. {\bf 131} (2003), no.4, 1151--1158.
\bibitem{Le1}{\sc J. Lehrb\"ack},
{\it Pointwise Hardy inequalities and uniformly fat sets},
Proc. Amer. Math. Soc. {\bf 136} (2008), no.6, 2193--2200.
\bibitem{L}{\sc J.L. Lewis},
{\it Uniformly fat sets},
Trans. Amer. Math. Soc. {\bf 308} (1988), no.1, 177--196.
\bibitem{M}{\sc V.G. Maz'ya}, 
``Sobolev spaces'',
Springer-Verlag, Berlin, 1985.
\bibitem{Ne}{\sc J. Ne\v{c}as},
{\it Sur une m\'ethode pour r\'esoudre les \'equations aux d\'eriv\'ees
partielles du type elliptique, voisine de la variationnelle}, 
Ann. Scuola Norm. Sup. Pisa (3) {\bf 16} (1962), 305--326.
\bibitem{Sh2}{\sc N. Shanmugalingam},
{\it Newtonian spaces: An extension of Sobolev spaces to metric
  measure spaces},
Rev. Mat. Iberoamericana, {\bf 16} (2000), 243--279.
\bibitem{Sh3}{\sc N. Shanmugalingam} 
{\it Harmonic functions on metric spaces},
Illinois J. Math. {\bf 45} (2001), no.3, 1021--1050. 
\bibitem{W}{\sc A. Wannebo},
{\it Hardy inequalities},  
Proc. Amer. Math. Soc. {\bf 109} (1990), no.1, 85--95.
\end{thebibliography}
\end{document}